%% file: regpol.tex
\documentclass[10pt,letterpaper]{article}
\usepackage{enumitem}
\usepackage{hyperref}
\usepackage{longtable}
\usepackage{amssymb,amsmath,amsfonts}
\usepackage{graphicx}
\graphicspath{{/}}
\usepackage[caption=false]{subfig}
\usepackage{listings}
\usepackage{anysize}
\marginsize{1.0in}{1.0in}{1.0in}{1.0in} 
\usepackage{booktabs}
\usepackage{url}
\usepackage{caption}
\newcommand{\besselroot}{j_\text{\tiny 0\hspace{-0.2ex}1}} \newcommand{\lambdaodot}{\besselroot^2}
\newcommand{\lambdaodotsquared}{\besselroot^4}
\newcommand{\lambdaodotcubed}{\besselroot^6}
\newcommand\maxS{{\mbox{150}}}

\newcommand{\GSB}{{\setlength\fboxsep{1pt}\mbox{GSB}}}
\def\ts[#1][#2]{\widetilde{{\{\mathbf{s}_{#1}}\}}_{#2}}
\def\s[#1][#2]{{\{\mathbf{s}_{#1}}\}_{#2}}
\def\tx[#1]{\widetilde{\mathbf{x}}_{#1}}
\def\ty[#1]{\widetilde{\mathbf{y}}_{#1}}
\def\tz[#1]{\widetilde{\mathbf{z}}_{#1}} \def\r[#1]{\mathbf{r}_{#1}}

\def\lambdaprime{\widehat\lambda}

\def\x[#1]{\mathbf{x}_{#1}}
\def\y[#1]{\mathbf{y}_{#1}}
\def\z[#1]{\mathbf{z}_{#1}}

\newcommand\markhere{\raisebox{-0.75ex}{\makebox[0pt]{\tiny$\blacktriangle$}}}

\begin{document}
\raggedright
\title{The fundamental Laplacian eigenvalue of the\\ regular polygon with Dirichlet boundary conditions}
\author{Robert Stephen Jones\thanks{\href{mailto:rsjones7@yahoo.com}{rsjones7@yahoo.com} , %
    \url{https://www.hbeLabs.com}}\\[0ex]%
  {\small Independent Researcher, Sunbury, Ohio}}
\maketitle

\begin{abstract}
The lowest eigenvalue of the Laplacian within the \mbox{$S$-sided}
regular polygon with Dirichlet boundary conditions is the focus of
this report.  As suggested by others, this eigenvalue may be expressed
as an asymptotic expansion in powers of~$1/S$ where,
interestingly, they have shown that the first few coefficients in that
expansion, up to sixth order, may be expressed analytically in terms
of Riemann zeta functions and roots of Bessel functions. This report
builds on that work with three main contributions: (1)~compelling
numerical evidence independently supporting those published results,
(2)~a conjecture adding two more terms to the asymptotic expansion,
and (3)~an observation that higher-order coefficients both alternate
in sign and grow rapidly in magnitude, which suggest the series
doesn't converge unless \mbox{$S\!\ge\!10$}.  This report is based on
a numerical computation of the eigenvalues precise to fifty digits
for~$S$ up to~\maxS.

\vspace{2ex}
\end{abstract}

\section*{Introduction}
Let $\Omega_S$ denote the interior of an $S$-sided regular polygon,
and $\partial\Omega_S$ its boundary.  The Laplacian eigenvalue problem
with Dirichlet boundary conditions for that polygon is defined by
\begin{equation}
  \left. \begin{array}{rl}
   \left[ \Delta + \lambda \right] \Psi(\r[]) = 0\quad & \text { for
     $\r[]\in\Omega_S$} \\[1ex]
     \Psi(\r[]) = 0 \quad & \text{ for
     $\r[]\in\partial\Omega_S$} 
  \end{array}\right\}
  \label{eq:problemstatement}
\end{equation} 
where $\Delta$ is the two-dimensional Laplacian, $\lambda>0$ is an
eigenvalue, and $\Psi(\r[])\not\equiv 0$ is a corresponding eigenfunction. 

A given regular polygon has an infinite tower of eigenvalues, but this
report shall focus only on the lowest (fundamental) eigenvalue, and,
more specifically, an asymptotic expansion of the form
\begin{equation}
  \lambda(S) \sim \lambdaprime(S)\equiv \lambdaodot\left[1+ \sum_{\mu=1}^{\infty} \frac{C_\mu}{S^\mu}\right] =
  \lambdaodot\left[1+\frac{C_1}{S}+\frac{C_2}{S^2}+\frac{C_3}{S^3}+\cdots\right]
\label{eq:series}
\end{equation}
where $\lambdaodot \approx 5.7831$ is the lowest
eigenvalue of the unit-radius circle\footnote{The corresponding
  [un-normalized] eigenfunction within that circle is
  $\Psi(\mathbf{r})=J_0(\besselroot r)$ where
  $r=|\mathbf{r}|\le 1$, and the number $\besselroot\approx
  2.4048$ is the first root of the Bessel function of the first
  kind,~$J_0(x)$.}.  Since we don't [yet] know how the
expansion converges, make the distinction: $\lambda(S)$ denotes the exact fundamental eigenvalue for
all $S=3,4,5,...$; whereas $\lambdaprime(S)$ denotes its asymptotic
expansion. When $\lambdaprime(S)$ is truncated at $C_N/S^N$, i.e., to $N$th
order, it shall be written $\lambdaprime^{[N]}(S)$ so that
\begin{equation}
       \lambda(S)\sim\lambdaprime(S)=\lambdaprime^{[N]}(S)+O(1/S^{N+1}).
\end{equation}

Built into the series is the assumption that as the polygon approaches
the unit-radius circle, i.e., in the limit $S\to\infty$,
\mbox{$\lambda(\infty)\!=\!\lambdaprime(\infty)\!=\!\lambdaodot$}. This
assumption seems natural, and is supported by numerics, but in
passing, recall the curious \emph{polygon-circle paradox} of
thin-plate theory~\cite{murray_1973} where an analogous assumption
breaks down.

When the area of $\Omega_S$ is held constant at $\pi$, i.e., the same
area as the unit-radius circle, the proposed expansion to eighth order is
\begin{multline}
  \lambdaprime^{[8]}(S)  = \lambdaodot\left\{ 1+
     \frac{4\,\zeta(3)}{S^3}
    + \frac{\left[ 12 - 2 \lambdaodot\right]\,\zeta(5)}{S^5}
    + \frac{\left[ 8 + 4 \lambdaodot\right]\,\zeta^2(3)}{S^6}\right. \\[1ex] \left.
             \mbox{} + \frac{\left[ 36-12\lambdaodot- \frac12\,\lambdaodotsquared\right]\,\zeta(7)}{S^7}
           + \frac{\left[ 48+8\lambdaodot+2\lambdaodotsquared\right]\,\zeta(3)\zeta(5)}{S^8}
            \right\} 
\label{eq:seriesfull}
\end{multline}
where $\zeta(n)=\sum_{\mu=1}^\infty\mu^{-n}$ is the well-known Riemann
zeta function, and where the last two terms are contributions of this
work.

The constant area of $\pi$ is chosen to more readily expose interesting facts about
the eigenvalue and its asymptotic expansion. Doing so automatically
factors out the well-known area rescaling dependence\footnote{If the
  polygon area is rescaled from $A$ to $tA$, an eigenvalue changes
  from $\lambda$ to $\lambda/t$, i.e., $A\lambda$ is constant.}  and
simplifies the expressions.  The Appendix details the relationship
between this \emph{transcribed} (equal area) regular polygon eigenvalue and the
inscribed one, as well as some other relationships.

Over the last twenty years, a few others have considered this
problem. A common theme is that those workers computed the
eigen-solution while gradually deforming the circle into the regular
polygon.  In 1997, Molinari \cite{m1997} suggested an expansion of
$\sqrt{\lambda(S)}$ in powers of $1/S$ and used conformal mapping to
estimate the leading coefficients of what he called a ``partial
resummation of terms in the $1/S$ expansion'', which he claimed
improved convergence.  In 2004, Grinfeld and Strang \cite{gs2004}
proposed\footnote{They did note, but without citation, that ``others''
  had already established that series.} Eq.~(\ref{eq:series}), and
they used ``the calculus of moving surfaces'' (CMS) and numerics to
estimate the first few coefficients.  Although of limited numerical
precision, these early efforts seemed promising and offered
interesting insights.

More recently, in 2012, Grinfeld and Strang \cite{gs2012} revisited
the problem and were able to express the coefficients up to $C_4$ as
integer multiples of the Riemann zeta function\footnote{If they
  actually used the constant $\pi$-area $\Omega_S$, which they
  suggested, they would have found
  $\lambdaprime^{[4]}(S)=\lambdaodot\{1+4\,\zeta(3)/S^3\}$. See the
  Appendix for details.}.  An interesting application of that work
appears in 2010 when
Oikonomou~\cite{o2010} studied the Casimir energy of a scalar field
within a regular polygon.  Several years later, in 2015,
Boady~\cite{phdthesisBoady}, working with Grinfeld, and also using
CMS, contributed two more terms, $C_5$ and~$C_6$. Their results were
obtained using a computer algebra system and do not depend on
numerical computations per~se. The terms up to sixth order~--~first
line of Eq.~(\ref{eq:seriesfull})~--~shall be referred to as the
Grinfeld-Strang-Boady [\GSB] terms.

Of note is that only two solutions
with finite $S$ are known in closed form, which for $\pi$-area $\Omega_S$ are
\begin{equation}
 \lambda(3) = \frac{4\pi}{\sqrt{3}} \qquad\text{and}\qquad \lambda(4)=2\pi
\end{equation}
All numerical evidence indicates that the $\pi$-area regular polygon eigenvalues
are monotonic with $S$,
\begin{equation}
  7.2552 \approx \frac{4\pi}{\sqrt{3}} \ge \lambda(S) > \lambda(S+1) > \lambdaodot\approx 5.7831
  \end{equation}
which is not unexpected~\cite{af2006,n2014}.

As terms are added to the asymptotic expansion per
Eq.~(\ref{eq:seriesfull}), interesting facts begin to emerge. For
example, the Riemann zeta function arguments are (so far) chosen from
\mbox{$\{3, 5, 7, ...\}$}, and~--~within each term~--~sum to that term's
order.\footnote{The next term in that sequence is most likely greater
  than eight but is otherwise not known.  Boady conjectured that the
  sequence consists of positive odd integers excluding~1, but the
  three numbers \mbox{$\{3, 5, 7\}$} also start the sequence of odd primes.}  That pattern
automatically requires \mbox{$C_1=C_2=C_4=0$}, which is a priori not obvious;
and, for example, the eighth order term involves only
$\zeta(3)\,\zeta(5)$ since $3+5$ is the only way to get~$8$ from that
set.  Of course, each added term also brings us a little closer to
identifying the elusive form of the function $\lambda(N)$ for which
$\lambdaprime(N)$ is merely its asymptotic expansion.

Of more practical interest is the ability to rapidly compute relatively
high-precision eigenvalues.
Indeed,
using Eq.~(\ref{eq:seriesfull}) and the computed, fifty-digit eigenvalues,
the relative discrepancy is empirically determined to be
\begin{equation}
  \frac{\lambdaprime^{[8]}(S) -
    \lambda(S)}{\lambda(S)}
        \approx \frac{-18.38}{S^{7.86}}
\label{eq:error}
\end{equation}
apparently valid for all $S\!\ge\!5$.  To illustrate, the ordinarily
difficult-to-calculate $S\!=\!128$ eigenvalue is readily found to a
precision of about fifteen digits,
\begin{align}
  \lambda(128) &=   5.78319922243209895\cdots \quad\text{(exact)}\\[-.5ex]
  \lambdaprime^{[8]}(128) &= 5.78319922243209606\cdots 
\end{align}
with a relative discrepancy of $-5.0\times 10^{-16}$, and where an
ellipsis in a number indicates truncation, not rounding.

My approach is quite straightforward. It begins with a high-precision
numerical computation~\cite{j2017} of the eigenvalues for~$S$ from~5
to~\maxS, precise to about fifty digits. These computed eigenvalues
(skipping the lowest few) are then fit using linear regression to a
truncated version of Eq.~(\ref{eq:series}) with just under forty
terms.  My conjecture for the seventh and eighth order terms is
derived using an LLL integer relation algorithm on the fit values of
the coefficients. The list of computed eigenvalues, the LLL technique,  
and some other details are provided in the Appendix.

All computations were performed on my personal commodity hardware
running free software\footnote{Typically per the GPL,
  \url{http://www.gnu.org/philosophy/free-sw.html}.} with GNU/Linux
(lubuntu 16.04.3~LTS) and its numerous ancillary utilities.  For the
eigenvalue computations, which took several months, I used a six-core
(12-thread) i7-5820K~@~3.30\,GHz with 64\,GB~RAM computer. Software of
choice was the {\tt pari/gp}~\cite{PARI293} calculator (compiled with {\tt gmp}
and {\tt pthread}).  A few symbolic computations were performed using
{\tt maxima}~\cite{maxima5.32}.

\section*{Linear Regression}

With the computed eigenvalues, linear regression shall be used to seek
numerical values of the coefficients $C_\mu$ of Eq.~(\ref{eq:series}).
Because this process uses up to around forty coefficients and requires
up to several dozen digits of precision, this unusual application of
linear regression requires some computational caveats.

The numerically computed eigenvalues shall be denoted
$\Lambda^{\text{[up]}}(S)$ and $\Lambda^{\text{[dn]}}(S)$ for the
upper and lower bounds, respectively, or generically, $\Lambda(S)$ (which can
refer to either bound or their mean).
The relative difference between the computed bounds satisfies
\begin{equation}
(0.114)\times 10^{-50} < \frac{\Lambda^\text{[up]}(S)-\Lambda^\text{[dn]}(S)}{\Lambda(S)} < (0.998)\times 10^{-50}
\end{equation}
which is just under~$10^{-50}$, and with a mean of~$(0.871)\times 10^{-50}$.

To develop the model equation, first let the independent variable be
$X=1/S$. The dependent variable $Y$ shall incorporate (1)~the computed
eigenvalues, (2)~the assumption that $\lambda(\infty)\!=\!\lambdaodot$,
and (3)~analytic expressions for the coefficients. Initially, all
coefficients are assumed unknown. Only after compelling numerical
evidence supports an analytic expression for a coefficient shall that
coefficient be considered known and exact, embellished with a tilde
(so that $C_\mu$ becomes $\widetilde{C}_\mu$), and incorporated
into~$Y$.

Virtually nothing is published regarding the convergence of
Eq.~(\ref{eq:series}) except for the vague but obvious notion that
convergence improves as $S$ increases. If we fit using low values of
$S$ for which the asymptotic series doesn't converge, the method will
fail because it won't capture the true nature of the function
$\lambda(S)$.  Therefore, the fit shall exclude the lowest few
computed eigenvalues, and that fit used to conjecture convergence
properties. 

To make this work, as many coefficients as possible must be included
in the fit, but not so many that the fit function begins to oscillate
wildly as it tries to ``connect the dots'' with a polynomial in~$X$.
Also, because of the numerically ill-conditioned nature of the linear
regression matrix computations, sufficient precision must be used.  To
that end, the precision of the linear regression computations is set to
a more-than-adequate 200 digits.  Both the number of terms to include in the
expansion and the computation precision are established experimentally.

\newcommand\myapprox{\approx}

To estimate the precision of the coefficients, the upper and lower
eigenvalue bounds are separately fit to the same model equation. This
process yields two sets of numerical values for the coefficients,
$\{C_\mu^\text{[up]}\}$ and $\{C_\mu^\text{[dn]}\}$, which
incidentally do not form bounds. The relative difference and
approximate number of digits in agreement between a pair of these
numbers are, respectively,
\begin{equation}
  \epsilon_\mu=\frac{C_\mu^\text{[up]}-C_\mu^\text{[dn]}}
          {\frac{1}{2} \left[ C_\mu^\text{[up]} + C_\mu^\text{[dn]}\right]} 
          \qquad\text{and}\qquad
          d_\mu = -\log_{10}\left|\epsilon_\mu\right|
\end{equation}
In this report, a given coefficient shall be reported
as the average value rounded one digit beyond a rounded $d_\mu$,
along with the value of $d_\mu$.  By example, if
$C_{28}^\mathrm{[up]}= 1.26128551\times 10^{16}$ and
$C_{28}^\mathrm{[dn]}= 1.26175766\times 10^{16}$, then
this coefficient is reported as
\begin{equation}
  C_{28}\ \{d_{28}\} = 1.262\times 10^{16}\ \{3.4\}
\end{equation}
where $\epsilon_{28}\!=\!-3.74\!\times\! 10^{-4}$.  Parameters are adjusted
so that \mbox{$d_\mu\!>\!1$} in every case. 

Three important observations regarding the numerical values of
the coefficients~--~looking down the series~--~include a drop in precision,
an alternation in sign, and a growth in magnitude. These observations
are quantified below.

In order to satisfy all the criteria, the following is chosen
\begin{equation}
  \text{\fbox{Fit Parameters: (a) include up to $C_{38}$ and (b) use $S=13, 14, \cdots, 150$}}
  \label{eq:fitparameters}
\end{equation}
In hindsight, this will ensure that the minimum $S$-value is not too small,
enough terms are included in the fit, and that every \mbox{$d_\mu\!>\!1$}.

There shall be four passes, of which the first three successfully
establish the analytic set of coefficients depicted in
Eq.~(\ref{eq:seriesfull}).  The final pass is used to estimate the
remaining coefficients. With each pass, the number of unknown
coefficients decreases, increasing their numerical precision slightly.

\subsection*{Pass 1}
To begin, assume all of the first 38 coefficients are unknown
and fit the computed eigenvalue data to the truncated series
\begin{equation}
  Y_a \equiv S\cdot
  \left[
    \frac{\Lambda(S)}{\lambdaodot}-1
    \right] = C_1 + C_2 X+ C_3 X^2 + \cdots + C_{38} X^{37}
  \label{eq:fitPass1}
\end{equation}
When this is done, the leading nine coefficients~--~listed in
Table~\ref{tab:pass1}~--~dramatically reveal that $C_1$, $C_2$, and
$C_4$ range from thirty to forty orders of magnitude smaller than the
nearby non-zero coefficients.  Indeed, to the precision of the
computation, they are effectively zero, which offers compelling
numerical evidence in support of the \GSB\  result that
\begin{equation}
  \fbox{$\widetilde{C}_1=  \widetilde{C}_2=  \widetilde{C}_4 =0$}
  \label{eq:zerocoeff}
\end{equation}
\begin{table}[!htb]
  \begin{center}
    \caption{First pass results. The first nine and the last of 38
      coefficients, assuming none are initially known per
      Eq.~(\ref{eq:fitPass1}).  Note the [near-]zero coefficients
      marked with arrows.
      \label{tab:pass1}
 }
    \vspace{1ex}
    \begin{tabular}{rl}\toprule
      $\mu$ & \multicolumn{1}{l}{\qquad $C_\mu\ \{d_\mu\}$}\\
      \midrule
1 & $+0.0000000000000000000000000000000000000019\ \{1.4\}$\makebox[0pt][l]{\ \ $\leftarrow$} \\
2 & $-0.0000000000000000000000000000000000039\ \{1.4\}$\makebox[0pt][l]{\ \ $\leftarrow$} \\
3 & $+4.8082276126383771415989526460458038\ \{34.4\}$ \\
4 & $-0.0000000000000000000000000000023\ \{1.3\}$\makebox[0pt][l]{\ \ $\leftarrow$} \\
5 & $+0.44964098545032430901630041787\ \{27.9\}$ \\
6 & $+44.98497175863112456004906931\ \{27.3\}$ \\
7 & $-50.539324388135164383037966\ \{24.9\}$ \\
8 & $+200.872237801870351587037\ \{23.1\}$ \\
9 & $-317.7704850739388022226\ \{21.1\}$ \\
\multicolumn{2}{l}{\qquad $\cdots$}\\
38 & $+2.53\times 10^{21}\ \{1.6\}$ \\
\bottomrule
    \end{tabular}
  \end{center}
\end{table}

\subsection*{Pass 2}

Next, incorporate Eq.~(\ref{eq:zerocoeff}) into the dependent variable and
refit the computed eigenvalue data to the model equation
\begin{align}
  Y_b & \equiv S^3\cdot
  \left[
    \frac{\Lambda(S)}{\lambdaodot}-1
    \right] = C_3 + C_5 X^2 + C_6 X^3 + \cdots + C_{38} X^{35}
  \label{eq:fitpass2}
\end{align}
which now includes 35 terms in the expansion on the right hand side.
The first three
coefficients of the fit are then compared to the non-zero \GSB\  terms%
\begin{align}
   \left.
  \begin{aligned}
         C_3\ \{d_3\} & =  4.808227612 6383771415 9895264604 579996267\ \{38.0\}\\
  4\,\zeta(3) & = 4.808227612 6383771415 9895264604 579996\markhere3059\cdots 
  \end{aligned} 
  \right\} & \quad 36\ \text{digits}\\[1ex]
   \left.
  \begin{aligned}
         C_5\ \{d_5\} & = 0.4496409854 5032430901 6300416830 27\ \{30.8\} \\
         (12-2\lambdaodot)\,\zeta(5) & =
                 0.4496409854 5032430901 63004168\markhere29 603 \cdots 
  \end{aligned} 
  \right\} & \quad 28\ \text{digits}\\[1ex]
   \left.
   \begin{aligned}
         C_6\ \{d_6\} & \myapprox 44.98497175 8631124560 0490696602 3\ \{29.8\} \\
         (8+4\lambdaodot)\,\zeta(3)^2 & =
                         44.98497175 8631124560 049069660\markhere 9 94 \cdots 
  \end{aligned} 
  \right\} & \quad 29\ \text{digits}
\end{align}
where the number of digits in agreement is indicated.
This result provides compelling numerical evidence
supporting the remaining \GSB\  terms,
\begin{equation}
  \fbox{$\widetilde{C}_3 = 4\,\zeta(3)$} \qquad
  \fbox{$\widetilde{C}_5 = (12-2\lambdaodot)\,\zeta(5)$} \qquad
  \fbox{$\widetilde{C}_6 = (8+4\lambdaodot)\,\zeta(3)^2$}
  \label{eq:gsbresult}
\end{equation}

\subsection*{Pass 3}

Next, incorporate the full \GSB\ result, Eqs.~(\ref{eq:zerocoeff})
and~(\ref{eq:gsbresult}), into the dependent variable and refit the
computed eigenvalue data to the model equation
\begin{align}
  Y_c & \equiv S^7\cdot
  \left[
    \frac{\Lambda(S)}{\lambdaodot}-
    \left(
    1
    +\frac{\widetilde{C}_3}{S^3}
    +\frac{\widetilde{C}_5}{S^5}
    +\frac{\widetilde{C}_6}{S^6}
    \right)
    \right] = C_7 + C_8 X + \cdots + C_{38} X^{31}
  \label{eq:fitPass3}
\end{align}
which now includes 32 terms in the expansion on the right hand side.
Comparing the resulting numerical coefficients $C_7$ and
$C_8$ to the proposed expressions yields
\begin{align}
  \left.
  \begin{aligned}
      C_7\ \{d_7\} & = -50.53932438813516438303806289079\ \{30.4\} \\
  \left(36-12\lambdaodot-\mbox{$\frac12$}\lambdaodotsquared\right)\,\zeta(7)&= 
                     -50.539324388135164383038062890\markhere 904\cdots
  \end{aligned} 
  \right\} & \quad 29\ \text{digits}\\[1ex]
  \left.
  \begin{aligned}
    C_8\ \{d_8\} & = +200.87223780187035158705886400\ \{27.8\}  \\
     \left(48+8\lambdaodot+2\lambdaodotsquared\right)\,\zeta(3)\,\zeta(5)&=
                     +200.872237801870351587058864\markhere 190\cdots
  \end{aligned} 
  \right\} & \quad 27\ \text{digits}
\end{align}
which provide compelling numerical evidence in support of my conjecture,
\begin{equation}
  \fbox{$\widetilde{C}_7 = \left(36-12\lambdaodot-\frac12\lambdaodotsquared\right)\,\zeta(7)$} \qquad
  \fbox{$\widetilde{C}_8 = \left(48+8\lambdaodot+2\lambdaodotsquared\right)\,\zeta(3)\zeta(5)$}
  \label{eq:jonesresult}
\end{equation}
To \emph{discover} the above relationships, the numerical coefficients, $C_7$ and
$C_8$, are input into an LLL integer relation algorithm using
\begin{equation}
\begin{aligned}
  0 &= a_7 C_7 + (b_7 + c_7 \lambdaodot+ d_7 \lambdaodotsquared)\,\zeta(7) \\[.5ex]
  0 &= a_8 C_8 + (b_8 + c_8 \lambdaodot+ d_8 \lambdaodotsquared)\,\zeta(3)\,\zeta(5)
\label{eq:LLLansatz}
\end{aligned}
\end{equation}
where the respective four integers $\{a, b, c, d\}$ are sought.  The
form of the integer relation is guided by the \GSB\ result. Details,
including a computer program, are given in the Appendix.

\subsection*{Pass 4}

Finally, incorporate Eqs.~(\ref{eq:zerocoeff}),~(\ref{eq:gsbresult}), and
(\ref{eq:jonesresult})
into the dependent variable and refit the computed eigenvalue data to
the model equation
\begin{multline}
  Y_d \equiv S^9\cdot
  \left[
    \frac{\Lambda(S)}{\lambdaodot}-
    \left(
    1
    +\frac{\widetilde{C}_3}{S^3}
    +\frac{\widetilde{C}_5}{S^5}
    +\frac{\widetilde{C}_6}{S^6}
    +\frac{\widetilde{C}_7}{S^7}
    +\frac{\widetilde{C}_8}{S^8}
    \right)
    \right] \\
  = C_9 + C_{10} X + C_{11} X^2 + \cdots + C_{38} X^{29}\qquad
  \label{eq:fitPass4}
\end{multline}
which now includes 30 terms in the expansion on the right hand side.
The complete results of this fit are shown in Table~\ref{tab:coefficientsPass04}. 
\begin{table}[htb!]
  \centering
  \caption{Fourth pass results listing all thirty coefficients of the model Eq.~(\ref{eq:fitPass4}).
  Note the sign alternation, growth in magnitude of $C_\mu$, and decrease in precision, $d_\mu$.}
    \label{tab:coefficientsPass04}
  \begin{tabular}{rl|rl}
    \toprule
    \multicolumn{1}{c}{$\mu$}&
    \multicolumn{1}{l|}{\hspace{3em}$C_\mu\ \{d_\mu\}$}&
    \multicolumn{1}{c}{$\mu$}&
    \multicolumn{1}{l}{\hspace{3em}$C_\mu\ \{d_\mu\}$}\\\midrule
 9 & $-317.77048507393880222654502267\ \{27.5\}$  & 24 & $+6.590391\times 10^{12}\ \{5.7\}$  \\
10 & $+1816.7620988762759616659826\ \{25.1\}$  & 25 & $-4.19643\times 10^{13}\ \{5.0\}$  \\
11 & $-6016.33571769034682922143\ \{22.8\}$  & 26 & $+2.7470\times 10^{14}\ \{4.4\}$  \\
12 & $+25200.97379293246467587\ \{20.9\}$  & 27 & $-1.8436\times 10^{15}\ \{3.9\}$  \\
13 & $-93352.057545638041207\ \{19.0\}$  & 28 & $+1.262\times 10^{16}\ \{3.4\}$  \\
14 & $+395412.696177504392\ \{17.2\}$  & 29 & $-8.702\times 10^{16}\ \{3.0\}$  \\
15 & $-1718008.2767654300\ \{15.6\}$  & 30 & $+5.934\times 10^{17}\ \{2.7\}$  \\
16 & $+7970543.96349877\ \{14.2\}$  & 31 & $-3.89\times 10^{18}\ \{2.5\}$  \\
17 & $-38310267.955146\ \{12.8\}$  & 32 & $+2.36\times 10^{19}\ \{2.3\}$  \\
18 & $+192454613.5202\ \{11.5\}$  & 33 & $-1.27\times 10^{20}\ \{2.1\}$  \\
19 & $-1004632656.0\ \{10.3\}$  & 34 & $+5.80\times 10^{20}\ \{2.0\}$  \\
20 & $+5.447327793\times 10^{9}\ \{9.2\}$  & 35 & $-2.12\times 10^{21}\ \{1.9\}$  \\
21 & $-3.05943716\times 10^{10}\ \{8.2\}$  & 36 & $+5.77\times 10^{21}\ \{1.8\}$  \\
22 & $+1.7770589\times 10^{11}\ \{7.3\}$  & 37 & $-1.03\times 10^{22}\ \{1.8\}$  \\
23 & $-1.065749\times 10^{12}\ \{6.5\}$  & 38 & $+8.96\times 10^{21}\ \{1.7\}$  \\
    \bottomrule
  \end{tabular}
\end{table}

\section*{Convergence}
\begin{figure}[hb!]
\centering\raggedright
\includegraphics[width=\textwidth]{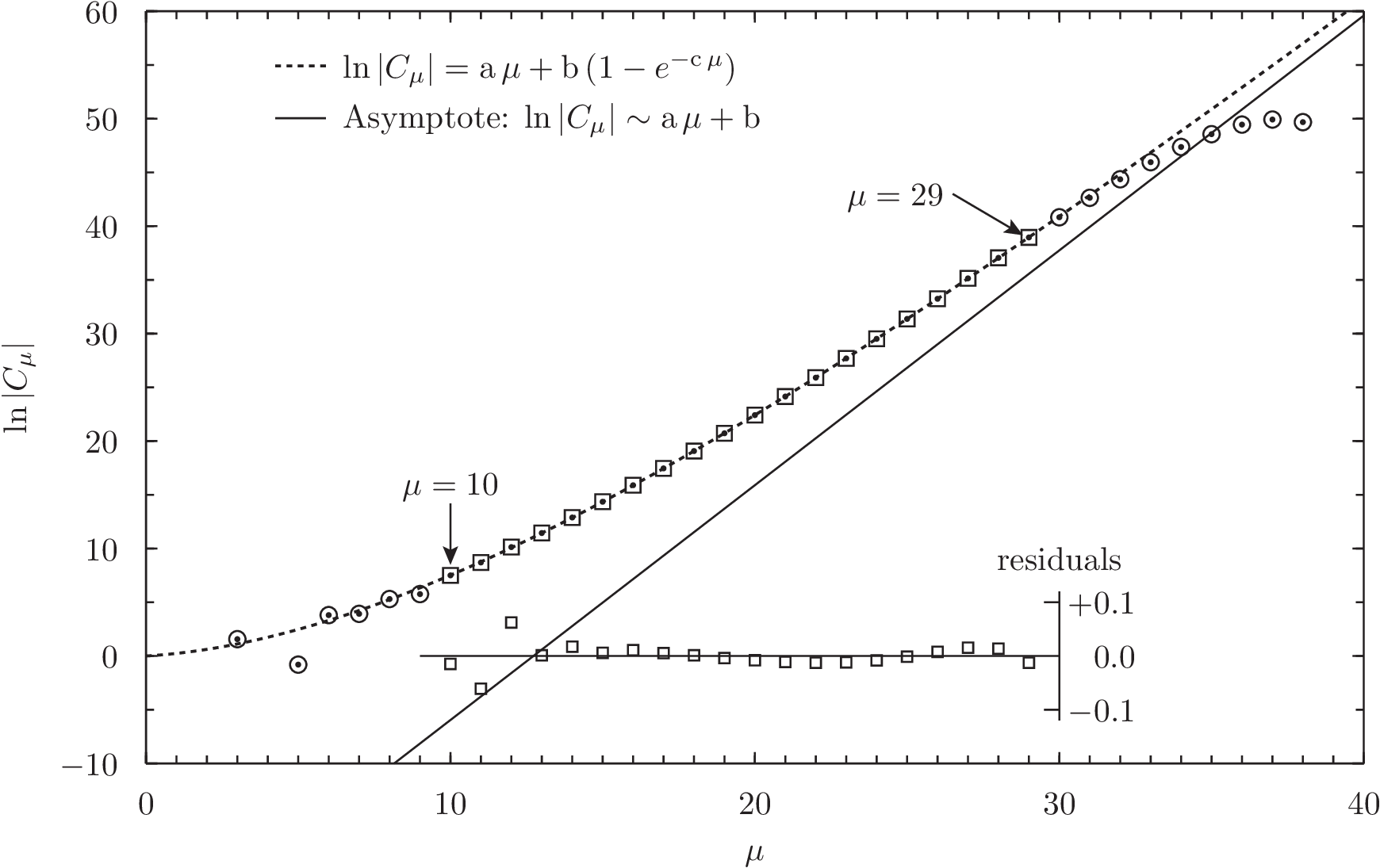}
\caption{Coefficient growth. The coefficient data (Pass~2) are fit to a simple
  model using from \mbox{$\mu=10$} to~$29$, inclusive. (Circled data
  excluded from fit.)  The inset displays the residuals of the points
  used in the fit on $50\times$
  scale.
\vspace{1ex}}
\label{fig.coeffgrowth}
\end{figure}
Numerically, the coefficients of Eq.~(\ref{eq:series}) 
exhibit two important properties which can be seen clearly in
Table~\ref{tab:coefficientsPass04}.
As one looks down that series, beyond the first few terms, the
coefficients appear
to both alternate in sign and grow in magnitude very rapidly,
apparently consistent with
\begin{align}
\begin{cases}
\displaystyle  C_\mu = (-1)^\mu |C_\mu|& \qquad\text{for $\mu>7$}\\
  \ln\left|C_\mu\right| \sim \mathrm{a}\,\mu+\mathrm{b}& \qquad\text{as $\mu\to\infty$}
\end{cases}
\end{align}
To analyze the coefficient growth, the natural logarithm of the
absolute value of the non-zero coefficients is plotted against the
series index in Fig.~\ref{fig.coeffgrowth} (circled and squared
dots). The down-turn in the last ten or so coefficients is
an artifact of the truncation (here, at $\mu\!=\!38$). Also visible (at
this scale) is that a smooth pattern isn't established until
$\mu\!\approx\!10$.

If ``$\ln\left|C_\mu\right|$'' does indeed approach a straight-line
asymptote, the slope~``$\mathrm{a}$'' of that asymptote determines
convergence. To see this, for sufficiently large but finite~$N$, the remainder
\begin{equation}
          R_N(S) \equiv \frac{\lambdaprime(S)-\lambdaprime^{[N\!-\!1]}(S)}{\lambdaodot}=\sum_{\mu=N}^{\infty}\frac{C_\mu}{S^\mu}
\approx  e^\mathrm{b} \sum_{i=N}^\infty\displaystyle (-1)^\mu \left[\frac{e^\mathrm{a}}{S}\right]^\mu
\end{equation}
certainly diverges if $S$ is too small. Indeed, the alternating geometric series in the last term
is absolutely convergent if $S > e^\mathrm{a}$, and diverges otherwise. (Both $\mathrm a$ and $S$ are positive.)
Since $\lambdaprime^{[N]}(S)$ is finite, conclude that
$\lambdaprime(S)$ is absolutely convergent if and only if
\begin{align}
  S > S_\text{cr}\equiv e^\mathrm{a} \qquad\text{(critical $S$-value)}
\end{align}
The challenge is to determine the slope of the asymptote, presuming
there is one. Without a known functional form for the coefficients,
there must be some mathematical caution with this numerical
exploration.  The most straightforward technique is to fit the
numerical values of the coefficients to a model with the simple
criteria that it have as few parameters as possible and a built-in
asymptote.

The model chosen here is the three-parameter exponential approach to
the asymptote,
\begin{equation}
  \ln\left|C_\mu\right| = \mathrm{a}\,\mu+\mathrm{b}\,\left(1-e^{-c\mu}\right)
\end{equation}
Some numerical values are excluded from the fit to avoid both the
artificial down-turn at the high end and the non-smooth behavior in
the low end. For the current set of data, using $10\!\le\!\mu\!\le\!29$ yields
\begin{equation}\left\{\begin{aligned}
               \mathrm{a} &= \phantom{+}2.185\pm 0.013\\ 
               \mathrm{b} &= -27.81\pm 0.57\\ 
               \mathrm{c} &= \phantom{+}0.07244\pm 0.0013
               \end{aligned}
\right.\end{equation}
where the expected values and standard deviations are reported. This
choice of model fits the data quite well and even extrapolates
through the lower end as shown in Fig~\ref{fig.coeffgrowth}. The inset
displays the residuals for the coefficients used in the fit.  With
these numbers, 
\begin{equation}
 S_\text{cr}= e^\mathrm{a} = 8.89\pm 0.12 = \underbrace{8.53\text{ to } 9.25}_{6 \sigma\ \text{interval}}
\end{equation}
Other simple models yield values anywhere between~7 and~9, but none as
large as~10. Erring on the side of caution, conclude that the
asymptotic series converges if $S$ is at least~$10$. In hindsight,
since $S_\text{cr}\!<\!13$, this result is consistent with the fit parameters,
Eq.~(\ref{eq:fitparameters}), used to determine the coefficients.

\section*{Future}

There is much room for future work.  For example, it is tempting to
search for yet higher-order coefficients and to study other
eigenvalues of the regular polygon. Another direction is to establish
more rigorous convergence criteria. Yet another higher goal is to
establish an analytic form of $\lambda(S)$, not merely its asymptotic
expansion.

Of note is that I am unable to extend the results to the ninth-order
term (or higher).  The natural extension of Eq.~\ref{eq:LLLansatz} might look like
\begin{equation}
 0 = a_9 C_9 + \left[b_9+c_9\lambdaodot+d_9\lambdaodotsquared+e_9\lambdaodotcubed\right]\,\zeta(3)^3
         +\left[f_9+g_9\lambdaodot+h_9\lambdaodotsquared+i_9\lambdaodotcubed\right]\,\zeta(9)
\end{equation}
where the nine [small] integers $\{a_9,b_9,...,i_9\}$ must be
determined~--~provided the Boady conjecture is somewhat valid. However, $C_9$ is
computed here to \emph{only} about 27 digits (Table~\ref{tab:coefficientsPass04}), and the LLL routine does
not suggest a unique solution as it does with the lower-order
terms. The failure may be due to either a breakdown in the simple
pattern or an insufficient precision for the LLL algorithm, or both.

\section*{Conclusion}
This investigation of the asymptotic expansion of the fundamental
Dirichlet eigenvalue of the Laplacian within the \mbox{$S$-sided}
regular polygon leads to three original results:
\begin{enumerate}[noitemsep,nolistsep]
\item independent and compelling numerical evidence in support of the \GSB\ result,
\item a conjecture for the next two terms (seventh and eighth order), and
\item numerical evidence that the asymptotic series may converge only
  if $S\!\ge\!10$.
\end{enumerate}   
These results are obtained using fifty-digit computed eigenvalues
for~$S$ up to~150.  Regression analysis of that data provides the
evidence in support of the \GSB\ result. The \GSB\ result, together
with an integer relation analysis of the numerical coefficients, leads
to the conjecture for the next two terms with compelling numerical
evidence supporting it.  Looking further down the asymptotic series, a
simple pattern (sign alternation and coefficient growth) emerged that
suggests it may converge only if \mbox{$S\!\ge\!10$}.

\section*{Appendix}
\textbf{Relation to the \GSB\ result}
The $S$-sided regular polygon used by
others~\cite{m1997,gs2004,o2010,gs2012,phdthesisBoady} is
typically \emph{inscribed} in a unit-radius circle. Grinfeld and
Strang use the term \emph{transcribe} to refer to an area-preserving
circle-to-polygon deformation.  Although I don't deform a circle, the
sequence ($S=3,4,5,...$) of $\pi$-area regular polygons shall herein
be referred to as \emph{transcribed} regular polygons to distinguish
them from \emph{inscribed} polygons, both in relation to that
unit-radius circle.

To distinguish the two problems, a prime is placed on the inscribed
problem variables.  Thus
\begin{equation}
 A(S)=\pi\qquad \text{ and }\qquad A'(S)= S \cos\left(\frac{\pi}{S}\right)\sin\left(\frac{\pi}{S}\right)
\end{equation}
are, respectively, the transcribed and inscribed area of the $S$-sided regular polygon.
Note that $A'(\infty)=\pi$, as required.

The area-rescaling relation for the eigenvalues is then either
\begin{equation}
  \lambda(S)\,\pi = \lambda'(S)A'(S) \qquad\text{ or }\qquad
  \lambdaprime(S)\,\pi = \lambdaprime'(S)\,A'(S)
\end{equation}
assuming $\lambdaprime(S)$ and $\lambdaprime'(S)$ converge. By
asymptotically expanding this, it is straightforward to show the
relationship between the published \GSB-result (inscribed) and my
$\pi$-area (transcribed) result, as discussed around Eqs.~11.46
to~11.48 of the Boady thesis.
The following, lightly-commented {\tt maxima}~\cite{maxima5.32} code will do that.
Note that my {\tt maxima}
function {\tt g(S)} is $\lambdaprime^{[6]}(S)/\lambdaodot$ and
variable {\tt L} is $\lambdaodot$.
\begin{lstlisting}
/* g(S) , for pi-area regular polygon, to sixth order */
g(S):=1 + 4*zeta(3)/S^3 + (12-2*L)*zeta(5)/S^5
                        + (8+4*L)*zeta(3)*zeta(3)/S^6;
/* Area of inscribed regular polygon divided by pi */
A(S):=(S/%pi)*cos(%pi/S)*sin(%pi/S);
/* expand out to sixth order in powers of 1/S */
expr:taylor(g(S)/A(S),S,inf,6);
/* Kludge to express in terms of zeta functions. */
expr:expr,[%pi^2=6*Z(2),%pi^4=90*Z(4),%pi^6=945*Z(6),
                 zeta(3)=Z(3),zeta(5)=Z(5)];
/* The next line shows the GSB result, Boady thesis Eq. 11.46 */
collectterms(expand(ratsimp(expr)),S);
\end{lstlisting}
The transcribed and inscribed expressions are, respectively,
\begin{align}
  \lambdaprime^{[6]}(S) & = \lambdaodot\,\left[ 1
  +\frac{4\,\zeta(3)}{S^3}
  + \frac{(12-2\lambdaodot)\,\zeta(5)}{S^5} 
   + \frac{(8+4\,\lambdaodot)\,\zeta^2(3)}{S^6} \right] \\
  \lambdaprime^{[6]\prime}(S) & = \lambdaodot \, \left[ 1
  +\frac{4\,\zeta(2)}{S^2} 
  +\frac{4\,\zeta(3)}{S^3}
  +\frac{28\,\zeta(4)}{S^4} \nonumber\right. \\
   + &\left. \frac{(12-2\lambdaodot)\,\zeta(5)+16\,\zeta(2)\,\zeta(3)}{S^5} 
   + \frac{(8+4\,\lambdaodot)\,\zeta^2(3)+124\,\zeta(6)}{S^6}\right]
\end{align}
One thing to note is that terms with even zeta function arguments
appear to be artifacts of the area dependence.

\textbf{Relation to other eigenvalues}

The focus of this report has been on the lowest eigenvalue.  It is
important to note that the asymptotic series may be readily modified
for eigenvalues within the same symmetry class as the lowest one,
i.e., those with $S$ even lines of symmetry intersecting at the
polygon center. Indeed, some of the other efforts were not limited to
the lowest eigenvalue. The extension is made by simply replacing
$\besselroot$ with the appropriate Bessel function root.  See,
for example, the Casimir energy analysis by
Oikonomou~\cite{o2010}.

\textbf{LLL procedure}

Below is a lightly-commented Pari gp-calculator program that applies
an LLL integer-relation algorithm to the five coefficients using the
gp routine {\tt qflll}.
\begin{lstlisting}
{
   default(realprecision,100); \\ working precision, 100 digits
   L0= (solve(x=sqrt(5.7),sqrt(5.8), besselj(0,x)))^2; \\ j_{01}^2
   o = [3,5,6,7,8]; \\ vector of term orders
   t = [1,2,2,3,3]; \\ number of terms in RHS

   C = vector(#o); \\ numerical coefficients          order
   C[1] = 4.80822761263837714159895264604579996267; \\ 3
   C[2] = 0.44964098545032430901630041683027;       \\ 5
   C[3] = 44.98497175863112456004906966023;         \\ 6
   C[4] = -50.53932438813516438303806289079;        \\ 7
   C[5] = 200.87223780187035158705886400;           \\ 8

   \\ Choosing zeta functions {3,5,6=3+3,7,8=3+5}
   Z = [zeta(3), zeta(5), zeta(3)^2, zeta(7), zeta(3)*zeta(5)];

   for (i=1, #o,  \\ loop over the five terms
      u = [C[i]]; \\ start off u-vector
      u = concat(u, Z[i]*vector(t[i], j, L0^(j-1)));
      N = #u;
      p = 30; \\ rounding parameter
      M = matid(N);  \\ NxN identity matrix
      M[N,]=round(u*10^p); \\ put numbers in last row
      v = qflll(M)[,1];    \\ LLL routine, 1st col returned
      v = sign(v[1])*v;    \\ make v[1] positive
      epsil = sum(j=1, #v, u[j]*v[j]); \\ ideally close to zero
      relerr= abs(epsil/C[i]); \\ relative error
      printf("C_%d=%6.3f ", o[i], C[i]);
      printf("relerr=%6.3g v=[%d", relerr, v[1]);
      for(j=1,t[i], printf(",%d", v[j+1])); printf("]\n");
   ); 
}
\end{lstlisting}
For reference, the output is
\begin{lstlisting}
C_3= 4.808 relerr=8.11 e-38 v=[1,-4]
C_5= 0.450 relerr=1.48 e-30 v=[1,-12,2]
C_6=44.985 relerr=1.70 e-29 v=[1,-8,-4]
C_7=-50.539 relerr=4.53 e-30 v=[2,-72,24,1]
C_8=200.872 relerr=9.46 e-28 v=[1,-48,-8,-2]
\end{lstlisting}
which can be used to construct the coefficients displayed in
Eq.~(\ref{eq:seriesfull}), including the non-zero \GSB\ terms.

\textbf{Computed eigenvalues}

The computed eigenvalues upon which this work relies are listed in
Table~\ref{tab:eigenvalues}.  This data represents a several-month
computation, from July~22 to November~5, 2017. The computer time
required for each eigenvalue increased with $S$ from a few seconds for
$S\!=\!5$ to about 2.5 days at~$S\!=\!150$. 

{
  \begin{longtable}{clc}
    \caption{Computed fifty-digit Dirichlet eigenvalues of the Laplacian,
      per Eq.~(\ref{eq:problemstatement}),
      for the $\pi$-area, $S$-sided
      regular polygon with~$S$ up to~$150$. By example,
      $\Lambda=5.12\,^{56}_{37}$ means $5.1237<\lambda<5.1256$.
      $\epsilon$ is the relative difference between upper and lower
      bounds. The first two eigenvalues are closed-form, and the last
      entry is the $S\!\to\!\infty$ circle eigenvalue.}
    \label{tab:eigenvalues}\\
    \toprule $S$ & \qquad $\Lambda$ & $\epsilon\ (\times 10^{-51})$ \\
    \midrule
    \endfirsthead
    \toprule $S$ & \qquad $\Lambda$ & $\epsilon\ (\times 10^{-51})$ \\
    \midrule 
    \endhead
    \multicolumn{3}{c}{$<<<$ Continued $>>>$}\\\bottomrule
    \endfoot
    {}\endlastfoot
    \input{eigenvalues.tex}
$\infty$ & \multicolumn{2}{l}{$5.783185962946784521175995758455807035071441806423685\cdots=\lambdaodot$}\\
    \bottomrule
\end{longtable}
}
\bibliographystyle{plain}
\bibliography{references}
\end{document}

%% file: eigenvalues.tex
$3$ & \multicolumn{2}{l}{$7.255197456936871402376313030568622929136264992370962\cdots = 4\pi/\sqrt{3}$} \\[1ex]
$4$ & \multicolumn{2}{l}{$6.283185307179586476925286766559005768394338798750211\cdots = 2\pi$} \\[1ex]
$5$ & $6.02213793204263387829800871005424296700530534044855\,^{88}_{18}$ & $1.14$ \\[1ex]
$6$ & $5.9174178316136612156885745768389615450082860040929\,^{721}_{266}$ & $7.67$ \\[1ex]
$7$ & $5.8664493126559858577124749417588410842427349136980\,^{702}_{495}$ & $3.52$ \\[1ex]
$8$ & $5.8384914335924428505166403795638157848367571520259\,^{684}_{419}$ & $4.52$ \\[1ex]
$9$ & $5.821826802270265731735546443716945921671786764620\,^{6082}_{5822}$ & $4.46$ \\[1ex]
$10$ & $5.8112603592191160227888164688111646234421581749002\,^{717}_{159}$ & $9.59$ \\[1ex]
$11$ & $5.8042306367174007218783944528561768184219904268906\,^{603}_{187}$ & $7.15$ \\[1ex]
$12$ & $5.7993698043565000793150253110077586111011868022191\,^{607}_{300}$ & $5.27$ \\[1ex]
$13$ & $5.7959002668560147097907710633371318342893868231820\,^{486}_{260}$ & $3.88$ \\[1ex]
$14$ & $5.793357005271194553273227078683828691973872008926\,^{8251}_{7934}$ & $5.45$ \\[1ex]
$15$ & $5.791450010651579975693848498149681163522889425220\,^{2343}_{1931}$ & $7.09$ \\[1ex]
$16$ & $5.7899918999902085343497522138280801741824641952139\,^{910}_{405}$ & $8.70$ \\[1ex]
$17$ & $5.7888578719811046986171966351899455269041934747495\,^{463}_{110}$ & $6.07$ \\[1ex]
$18$ & $5.7879625918578468642125683801538930922402544098186\,^{443}_{032}$ & $7.08$ \\[1ex]
$19$ & $5.787246351381961243008036644834744747679263989712\,^{6307}_{5843}$ & $8.01$ \\[1ex]
$20$ & $5.7866665141403722135309129620257370578364779208083\,^{802}_{290}$ & $8.84$ \\[1ex]
$21$ & $5.7861920775968442730282037573295124400580514989647\,^{594}_{040}$ & $9.56$ \\[1ex]
$22$ & $5.7858001294283650275745860443434083658465152450013\,^{521}_{126}$ & $6.81$ \\[1ex]
$23$ & $5.7854734864549016320482640701967564443197042731253\,^{729}_{441}$ & $4.96$ \\[1ex]
$24$ & $5.785199089790024091834463612979690633039598488051\,^{8322}_{7878}$ & $7.66$ \\[1ex]
$25$ & $5.7849668941304235014186706838851226120341859263223\,^{965}_{502}$ & $7.98$ \\[1ex]
$26$ & $5.7847690863148429779922747176926528381622936339967\,^{445}_{104}$ & $5.88$ \\[1ex]
$27$ & $5.7845995272364846405932228271557722785394731147071\,^{640}_{151}$ & $8.44$ \\[1ex]
$28$ & $5.7844533477517199511967948423990305225015801633632\,^{812}_{448}$ & $6.28$ \\[1ex]
$29$ & $5.784326652365411207380293385946299260682128920466\,^{5349}_{4846}$ & $8.68$ \\[1ex]
$30$ & $5.784216299392264044119036734072080273019850029152\,^{1071}_{0693}$ & $6.53$ \\[1ex]
$31$ & $5.784119736080032703344528384598811278461255613439\,^{3095}_{2588}$ & $8.75$ \\[1ex]
$32$ & $5.784034873702444318330507486930290862975888102019\,^{3253}_{2867}$ & $6.65$ \\[1ex]
$33$ & $5.783959992040508812335032522628621979446533803257\,^{2140}_{1637}$ & $8.69$ \\[1ex]
$34$ & $5.7838936656948092520334765687133149515574394886548\,^{910}_{410}$ & $8.62$ \\[1ex]
$35$ & $5.7838347067709882028407000052233665412473857264158\,^{612}_{227}$ & $6.64$ \\[1ex]
$36$ & $5.783782119955880627699919965851270952506815797405\,^{8288}_{7800}$ & $8.42$ \\[1ex]
$37$ & $5.783735067049846291962440636524259616676854699958\,^{3345}_{2864}$ & $8.30$ \\[1ex]
$38$ & $5.7836928387732922677065179222622535502538886661416\,^{934}_{461}$ & $8.16$ \\[1ex]
$39$ & $5.7836548322108711433862079109147809210199222742790\,^{551}_{087}$ & $8.02$ \\[1ex]
$40$ & $5.783620532655973576951368558889088191821701701593\,^{5398}_{4942}$ & $7.87$ \\[1ex]
$41$ & $5.783589498912728857243541087541613440608196223660\,^{2233}_{1681}$ & $9.53$ \\[1ex]
$42$ & $5.783561351331960679963744950469531484077839530288\,^{6533}_{5995}$ & $9.29$ \\[1ex]
$43$ & $5.7835357620219719712774467620660707772186262259777\,^{895}_{372}$ & $9.04$ \\[1ex]
$44$ & $5.783512446799268033474803357987667934572313002680\,^{7317}_{6808}$ & $8.79$ \\[1ex]
$45$ & $5.783491158538856302913858878804287379136291682036\,^{6143}_{5647}$ & $8.55$ \\[1ex]
$46$ & $5.7834716816561681101051372252375789382247863870892\,^{595}_{113}$ & $8.31$ \\[1ex]
$47$ & $5.783453827508463297379645913916576075154358289917\,^{3155}_{2592}$ & $9.72$ \\[1ex]
$48$ & $5.7834374305468654192506360647558848068434756772524\,^{645}_{100}$ & $9.41$ \\[1ex]
$49$ & $5.783422345083940177286945516586288853576803572305\,^{1365}_{0838}$ & $9.11$ \\[1ex]
$50$ & $5.7834084425682129947801161963236469590292923516017\,^{854}_{343}$ & $8.82$ \\[1ex]
$51$ & $5.7833956092779034421663981399445418601302161224952\,^{893}_{399}$ & $8.53$ \\[1ex]
$52$ & $5.7833837443627027000195590153219733045460685739016\,^{742}_{177}$ & $9.76$ \\[1ex]
$53$ & $5.7833727581755982440480496612883417072850762249903\,^{706}_{160}$ & $9.42$ \\[1ex]
$54$ & $5.7833625708472927428973141440030573413552170062\,^{900473}_{899947}$ & $9.09$ \\[1ex]
$55$ & $5.783353111064236385644810842163192250719443038200\,^{6093}_{5585}$ & $8.77$ \\[1ex]
$56$ & $5.78334431501812935463844781588756191043588875765\,^{40427}_{39854}$ & $9.89$ \\[1ex]
$57$ & $5.7833361255002921030903693546648597586536448335457\,^{618}_{066}$ & $9.53$ \\[1ex]
$58$ & $5.7833284911188091539136722748217426193090792309470\,^{845}_{313}$ & $9.18$ \\[1ex]
$59$ & $5.7833213656200338537379396114645869789760826156081\,^{794}_{282}$ & $8.84$ \\[1ex]
$60$ & $5.783314707299059428210797800183990949668498895119\,^{6345}_{5774}$ & $9.86$ \\[1ex]
$61$ & $5.783308478486244250571058736219661811034316213866\,^{3113}_{2564}$ & $9.48$ \\[1ex]
$62$ & $5.7833026450989284101703805978440170235282052458579\,^{531}_{002}$ & $9.13$ \\[1ex]
$63$ & $5.783297176249175699422050587208076227946408600731\,^{7365}_{6856}$ & $8.79$ \\[1ex]
$64$ & $5.7832920438997850274611258286193147488559177241149\,^{925}_{364}$ & $9.70$ \\[1ex]
$65$ & $5.783287222561990216101379318958387384571634250267\,^{3215}_{2675}$ & $9.32$ \\[1ex]
$66$ & $5.7832826890292492111478490456613641227808808137750\,^{732}_{212}$ & $8.97$ \\[1ex]
$67$ & $5.783278422142346980297970170232395388356378905112\,^{1470}_{0900}$ & $9.83$ \\[1ex]
$68$ & $5.783274402581728387648667336549538768975499369192\,^{6319}_{5771}$ & $9.45$ \\[1ex]
$69$ & $5.7832706126835606254665067931060003616486874965111\,^{554}_{028}$ & $9.08$ \\[1ex]
$70$ & $5.783267036276517720104953505291690754885783103053\,^{3438}_{2865}$ & $9.89$ \\[1ex]
$71$ & $5.7832636585366972675053570571575828141496986506519\,^{916}_{366}$ & $9.50$ \\[1ex]
$72$ & $5.7832604658584342703907945644030276005726248720419\,^{576}_{046}$ & $9.14$ \\[1ex]
$73$ & $5.7832574457390789460408703608938752178996635604860\,^{577}_{004}$ & $9.90$ \\[1ex]
$74$ & $5.7832545866760630843215845592704878980855033362418\,^{808}_{257}$ & $9.51$ \\[1ex]
$75$ & $5.7832518780747999519072844129240886957218156591922\,^{670}_{141}$ & $9.14$ \\[1ex]
$76$ & $5.7832493101661516716346291929211527006058362175622\,^{945}_{374}$ & $9.85$ \\[1ex]
$77$ & $5.7832468739323602985306185231341800220529994071590\,^{639}_{150}$ & $8.45$ \\[1ex]
$78$ & $5.783244561040478512305563118717910086346504177634\,^{6472}_{5945}$ & $9.09$ \\[1ex]
$79$ & $5.7832423637824563386561208792114744191451563621838\,^{703}_{137}$ & $9.76$ \\[1ex]
$80$ & $5.7832402750211444434602632903336057637956537641100\,^{733}_{246}$ & $8.41$ \\[1ex]
$81$ & $5.783238288141564708788818913744019150004231228987\,^{6027}_{5504}$ & $9.01$ \\[1ex]
$82$ & $5.7832363970068770166311668301521069133506020689353\,^{844}_{285}$ & $9.64$ \\[1ex]
$83$ & $5.78323459591853914193516915928294052374939682204\,^{30082}_{29599}$ & $8.34$ \\[1ex]
$84$ & $5.78323287958021583661957400684452483337094433777\,^{80111}_{79539}$ & $9.87$ \\[1ex]
$85$ & $5.783231243065044797761869159688729531153474161960\,^{8425}_{7929}$ & $8.57$ \\[1ex]
$86$ & $5.7832296817859122999872941895330211636707297244341\,^{803}_{275}$ & $9.12$ \\[1ex]
$87$ & $5.7832281914684307238009975309857996746289319161078\,^{918}_{356}$ & $9.69$ \\[1ex]
$88$ & $5.783226768126344787904121271120141989640291161915\,^{8089}_{7600}$ & $8.44$ \\[1ex]
$89$ & $5.783225408039123644386441987435525040163562180818\,^{3249}_{2730}$ & $8.96$ \\[1ex]
$90$ & $5.783224107731522678222767981643965415039391863402\,^{5443}_{4893}$ & $9.49$ \\[1ex]
$91$ & $5.783222863954922345132559440101223356983832893756\,^{4446}_{3917}$ & $9.12$ \\[1ex]
$92$ & $5.783221673670272096140477542110521995453407740228\,^{7388}_{6879}$ & $8.78$ \\[1ex]
$93$ & $5.7832205340324857277389860275897345239229247945415\,^{675}_{138}$ & $9.27$ \\[1ex]
$94$ & $5.7832194423761506695184080213013434106581037042669\,^{600}_{033}$ & $9.78$ \\[1ex]
$95$ & $5.7832183962024280411626163932933568556831624568013\,^{978}_{433}$ & $9.41$ \\[1ex]
$96$ & $5.783217393167033007010852241398377430861882282949\,^{2327}_{1753}$ & $9.90$ \\[1ex]
$97$ & $5.783216431069196227758518741530320098257066954853\,^{8270}_{7765}$ & $8.71$ \\[1ex]
$98$ & $5.783215507841517228010683796245184675675966031085\,^{8128}_{7597}$ & $9.16$ \\[1ex]
$99$ & $5.783214621540629415561788196159738903682817178900\,^{4118}_{3560}$ & $9.62$ \\[1ex]
$100$ & $5.783213770338604434367144876070106403216402193565\,^{4364}_{3827}$ & $9.26$ \\[1ex]
$101$ & $5.7832129525150306224265566440583980748230299905146\,^{534}_{018}$ & $8.91$ \\[1ex]
$102$ & $5.783212166449706678032327479127825448347817055273\,^{1180}_{0639}$ & $9.34$ \\[1ex]
$103$ & $5.783211410615897300382237922986464486484414973146\,^{3446}_{2880}$ & $9.78$ \\[1ex]
$104$ & $5.783210683574102639949628780569819360526210833755\,^{8062}_{7560}$ & $8.66$ \\[1ex]
$105$ & $5.783209983966297937380176592129591988164382058292\,^{3293}_{2723}$ & $9.84$ \\[1ex]
$106$ & $5.783209310510603806034215145692410082844657773366\,^{6300}_{5794}$ & $8.72$ \\[1ex]
$107$ & $5.7832086619963512745011282308896272842754873361699\,^{771}_{243}$ & $9.11$ \\[1ex]
$108$ & $5.7832080372795089971773313237835737211106455961164\,^{988}_{437}$ & $9.51$ \\[1ex]
$109$ & $5.7832074352784430036203824301621058164531928714376\,^{652}_{077}$ & $9.92$ \\[1ex]
$110$ & $5.783206854969982026448346069044505432273217666833\,^{1081}_{0569}$ & $8.82$ \\[1ex]
$111$ & $5.783206295385763854495868384691346982424634552152\,^{1382}_{0850}$ & $9.19$ \\[1ex]
$112$ & $5.783205755608840330598899278697864669617886870683\,^{4198}_{3644}$ & $9.57$ \\[1ex]
$113$ & $5.7832052347705205764174624901385415125593493624890\,^{856}_{280}$ & $9.95$ \\[1ex]
$114$ & $5.7832047320474338019919255227789909098567095792232\,^{862}_{347}$ & $8.88$ \\[1ex]
$115$ & $5.7832042466587946646835703161902916564434900898372\,^{578}_{002}$ & $9.95$ \\[1ex]
$116$ & $5.783203777863855598037977272725438352252698382214\,^{6403}_{5888}$ & $8.89$ \\[1ex]
$117$ & $5.78320332495953185129119383617831568697480390758\,^{30315}_{29780}$ & $9.23$ \\[1ex]
$118$ & $5.7832028872781861783972732175519624981091930701562\,^{574}_{019}$ & $9.58$ \\[1ex]
$119$ & $5.7832024641855612037927177369523705379144297168954\,^{592}_{017}$ & $9.93$ \\[1ex]
$120$ & $5.7832020550788484815190743450796360573620493874627\,^{710}_{156}$ & $9.57$ \\[1ex]
$121$ & $5.783201659384884164529487670253652645264573710311\,^{4258}_{3724}$ & $9.22$ \\[1ex]
$122$ & $5.7832012765584620207102095560653892586849534945666\,^{979}_{426}$ & $9.55$ \\[1ex]
$123$ & $5.783200906080755279143481813903136517206041981035\,^{3386}_{2814}$ & $9.87$ \\[1ex]
$124$ & $5.783200547457839471403393668939526347739671446805\,^{1004}_{0452}$ & $9.52$ \\[1ex]
$125$ & $5.7832002002193090544707614611862501293189148764360\,^{754}_{184}$ & $9.84$ \\[1ex]
$126$ & $5.783199863916981169795599727550027796605983950656\,^{5444}_{4894}$ & $9.48$ \\[1ex]
$127$ & $5.783199538123680412174552013868591251347340966166\,^{6035}_{5505}$ & $9.15$ \\[1ex]
$128$ & $5.783199222432098956985238320133014718214303003483\,^{5340}_{4793}$ & $9.44$ \\[1ex]
$129$ & $5.783198916453726829015452454213861349980377100531\,^{3297}_{2733}$ & $9.74$ \\[1ex]
$130$ & $5.783198619817847494322697718282907407501883505901\,^{3379}_{2835}$ & $9.40$ \\[1ex]
$131$ & $5.783198332170594321576460839361900033150683319289\,^{3283}_{2758}$ & $9.07$ \\[1ex]
$132$ & $5.783198053174063794165401631077445539852104232715\,^{8325}_{7747}$ & $9.98$ \\[1ex]
$133$ & $5.7831977825054816616973973974236090228425636259915\,^{783}_{225}$ & $9.63$ \\[1ex]
$134$ & $5.783197519856418501824964777729256102794211740624\,^{6034}_{5460}$ & $9.91$ \\[1ex]
$135$ & $5.78319726493205142280133747938195184354472751162\,^{40300}_{39781}$ & $8.97$ \\[1ex]
$136$ & $5.783197017450468875813125929970793135181216310333\,^{2208}_{1673}$ & $9.23$ \\[1ex]
$137$ & $5.7831967771420157657553791022485248591356074629077\,^{779}_{228}$ & $9.50$ \\[1ex]
$138$ & $5.783196543748676251355969965179691633199779643683\,^{3230}_{2699}$ & $9.17$ \\[1ex]
$139$ & $5.78319631702349181190803820593767295255592560127\,^{10466}_{09919}$ & $9.43$ \\[1ex]
$140$ & $5.7831960967300123296839310805154072611401163978360\,^{776}_{215}$ & $9.69$ \\[1ex]
$141$ & $5.7831958826417780956100870897864370536602334172115\,^{593}_{016}$ & $9.95$ \\[1ex]
$142$ & $5.78319567454183079209634447072860419375002513070\,^{20001}_{19477}$ & $9.03$ \\[1ex]
$143$ & $5.7831954722222516420511547774247258923709033792146\,^{705}_{134}$ & $9.86$ \\[1ex]
$144$ & $5.783195275483725038000991945505785473982967396266\,^{5381}_{4829}$ & $9.52$ \\[1ex]
$145$ & $5.783195084135126080710376077178119900703268965184\,^{8025}_{7492}$ & $9.20$ \\[1ex]
$146$ & $5.783194897993130563536383404875214594259743023232\,^{7540}_{6993}$ & $9.43$ \\[1ex]
$147$ & $5.7831947168818460376488389106165087799453564413638\,^{878}_{318}$ & $9.67$ \\[1ex]
$148$ & $5.78319454063246268484401771071146302825917616988\,^{20197}_{19656}$ & $9.34$ \\[1ex]
$149$ & $5.7831943690829228095596336817453011340629123532482\,^{985}_{431}$ & $9.57$ \\[1ex]
$150$ & $5.7831942020776078403959106662108900488674959487428\,^{954}_{386}$ & $9.80$ \\[1ex]